\documentclass[a4paper,11pt]{article}
\usepackage{latexsym}
\usepackage{amssymb}
\usepackage{enumerate}
\usepackage{amsfonts}
\usepackage{amsmath}
\usepackage[dvips]{graphicx}
\usepackage[paper=a4paper,left=30mm,right=20mm,top=25mm,bottom=30mm]{geometry}
    \newenvironment{dedication}
        {\vspace{2ex}\begin{quotation}\begin{center}\begin{em}}
        {\par\end{em}\end{center}\end{quotation}}
\newcommand{\qed}{$\Box$}

\newenvironment{@abssec}[1]{%
    \if@twocolumn

      \section*{#1}%
    \else

      \vspace{.05in}\footnotesize
      \parindent .2in
 {\upshape\bfseries #1. }\ignorespaces
    \fi}

    {\if@twocolumn\else\par\vspace{.1in}\fi}

\newenvironment{keywords}{\begin{@abssec}{\keywordsname}}{\end{@abssec}}

\newenvironment{AMS}{\begin{@abssec}{\AMSname}}{\end{@abssec}}

\newcommand\keywordsname{Key words}
\newcommand\AMSname{AMS subject classifications}
\newcommand\AMname{AMS subject classification}
\newtheorem{theorem}{Theorem}
 \newtheorem{lemma}[theorem]{Lemma}
 \newtheorem{corollary}[theorem]{Corollary}
 \newtheorem{proposition}[theorem]{Proposition}

\def\qed{\vbox{\hrule height0.6pt\hbox{%
  \vrule height1.3ex width0.6pt\hskip0.8ex
  \vrule width0.6pt}\hrule height0.6pt
 }}

\title{Two-phase heat conductors with a stationary isothermic surface
\thanks{This research was partially supported by the Grant-in-Aid
for Scientific Research (B) ($\sharp$ 26287020) of
Japan Society for the Promotion of Science.}}

\author{Shigeru Sakaguchi\thanks{Research Center for Pure and Applied Mathematics,
Graduate School of  Information Sciences, Tohoku
University, Sendai, 980-8579,  Japan.
({\tt sigersak@m.tohoku.ac.jp}).}}
\begin{document}

\maketitle

\begin{dedication}
\begin{large}
Dedicated to Professor Giovanni Alessandrini on his sixties birthday
\end{large}
\end{dedication}

\begin{abstract}
We consider a two-phase heat conductor in $\mathbb R^N$ with $N \geq 2$ consisting of a core and a shell with different constant conductivities. Suppose that, initially, the conductor has temperature 0 and, at all times, its boundary is kept at temperature 1.
 It is shown that, if  there is a stationary isothermic surface in the shell near the boundary, then the structure of the conductor must be spherical. Also, when the medium outside the two-phase conductor has a possibly different conductivity,  we consider the Cauchy problem with $N \ge 3$ and  the initial condition where the conductor has temperature 0 and  the outside medium has temperature 1. Then we show that almost the same proposition holds true.
 \end{abstract}




\pagestyle{plain}
\thispagestyle{plain}

\section{Introduction}
\label{introduction}

\vskip 2ex
Let $\Omega$ be a bounded $C^2$ domain in $\mathbb R^N\ (N \ge 2)$ with boundary $\partial\Omega$, and let $D$ be a bounded  $C^2$ open set in $\mathbb R^N$ which may have finitely many connected components.  Assume that $\Omega\setminus\overline{D}$ is connected and $\overline{D} \subset \Omega$. Denote by $\sigma=\sigma(x)\ (x \in \mathbb R^N)$  the conductivity distribution of the medium given by
$$
\sigma =
\begin{cases}
\sigma_c \quad&\mbox{in } D, \\
\sigma_s \quad&\mbox{in } \Omega \setminus D, \\
\sigma_m \quad &\mbox{in } \mathbb R^N \setminus \Omega,
\end{cases}
$$
where $\sigma_c, \sigma_s, \sigma_m$ are positive constants and $\sigma_c \not=\sigma_s$. This kind of three-phase electrical conductor has been dealt with in \cite{KLS2016} in the study of neutrally coated inclusions.

In the present paper we consider the heat diffusion over two-phase or three-phase heat conductors.
Let $u =u(x,t)$ be the unique bounded solution of either the initial-boundary value problem for the diffusion equation:
\begin{eqnarray}
&&u_t =\mbox{ div}( \sigma \nabla u)\quad\mbox{ in }\ \Omega\times (0,+\infty), \label{heat equation initial-boundary}
\\
&&u=1  \ \quad\qquad\qquad\mbox{ on } \partial\Omega\times (0,+\infty), \label{heat Dirichlet}
\\ 
&&u=0  \  \quad\qquad\qquad \mbox{ on } \Omega\times \{0\},\label{heat initial}
\end{eqnarray}
or the Cauchy problem for the diffusion equation:
\begin{equation}
  u_t =\mbox{ div}(\sigma \nabla u)\quad\mbox{ in }\  \mathbb R^N\times (0,+\infty) \ \mbox{ and }\ u\ ={\mathcal X}_{\Omega^c}\ \mbox{ on } \mathbb R^N\times
\{0\},\label{heat Cauchy}
\end{equation}
where ${\mathcal X}_{\Omega^c}$ denotes the characteristic function of the set $\Omega^c=\mathbb R^N \setminus\Omega$. Consider a bounded domain $G$ in $\mathbb R^N$ satisfying
\begin{equation}
\label{near the boundary}
\overline{D} \subset G \subset \overline{G} \subset \Omega\ \mbox{ and } \mbox{ dist}(x,\partial\Omega) \le \mbox{ dist}(x, \overline{D})\ \mbox{ for every } x \in \partial G.
\end{equation}
The purpose of the present paper is to show the following theorems.
\begin{theorem} 
\label{th:stationary isothermic} Let $u$ be the solution of problem \eqref{heat equation initial-boundary}-\eqref{heat initial} for $N \ge 2$, and let
$\Gamma$ be a connected component of $\partial G$ satisfying
\begin{equation}
\label{nearest component}
\mbox{\rm dist}(\Gamma, \partial\Omega) = \mbox{\rm dist}(\partial G, \partial\Omega).
\end{equation}  
If there exists a function $a : (0, +\infty) \to (0, +\infty) $ satisfying
\begin{equation}
\label{stationary isothermic surface partially}
u(x,t) = a(t)\ \mbox{ for every } (x,t) \in \Gamma \times (0, +\infty),
\end{equation}
then $\Omega$ and $D$ must be concentric balls.
\end{theorem}

\begin{corollary} 
\label{cor:stationary isothermic} Let $u$ be the solution of problem \eqref{heat equation initial-boundary}-\eqref{heat initial} for $N \ge 2$.  If
there exists a function $a : (0, +\infty) \to (0, +\infty) $ satisfying
\begin{equation}
\label{stationary isothermic surface}
u(x,t) = a(t)\ \mbox{ for every } (x,t) \in \partial G\times (0, +\infty),
\end{equation}
then $\Omega$ and $D$ must be concentric balls.
\end{corollary}

\begin{theorem} 
\label{th:stationary isothermic cauchy} Let $u$ be the solution of problem \eqref{heat Cauchy} for $N \ge 3$. Then the following assertions hold:
\begin{itemize}
\item[\rm (1)] If
there exists a function $a : (0, +\infty) \to (0, +\infty) $ satisfying \eqref{stationary isothermic surface},
then $\Omega$ and $D$ must be concentric balls.
\item[\rm (2)] If $\sigma_s=\sigma_m$ and there exists a function $a : (0, +\infty) \to (0, +\infty) $ satisfying \eqref{stationary isothermic surface partially} for a connected component $\Gamma$ of $\partial G$ with \eqref{nearest component}, then  $\Omega$ and $D$ must be concentric balls.
\end{itemize}
\end{theorem}

Corollary \ref{cor:stationary isothermic} is just an easy by-product of Theorem \ref{th:stationary isothermic}. Theorem \ref{th:stationary isothermic cauchy} is limited to the case where $N \ge 3$, which is not natural; that is required for  technical reasons in the use of the auxiliary functions $U, V, W$ given in section \ref{section4}. We conjecture that Theorem \ref{th:stationary isothermic cauchy} holds true also for $N=2$.

The condition \eqref{stationary isothermic surface partially} means that $\Gamma$ is an isothermic surface  of the normalized temperature $u$ at every time, and hence $\Gamma$ is called a {\it stationary} isothermic surface of $u$.
When $D = \emptyset$ and $\sigma$ is constant on $\mathbb R^N$,  a symmetry theorem similar to Theorem \ref{th:stationary isothermic} or Theorem \ref{th:stationary isothermic cauchy} has been proved in \cite[Theorem 1.2, p. 2024]{MSmmas2013} provided the conclusion is replaced by that $\partial\Omega$  must be either a sphere or the union of two concentric spheres, and a symmetry theorem similar to Corollary \ref{cor:stationary isothermic} has also  been proved in \cite[Theorem 1.1, p. 932]{MSannals2002}. The present paper gives a generalization of the previous results to multi-phase heat conductors. 

We note that the study of the relationship between the stationary isothermic surfaces and the symmetry of the problems has been initiated by Alessandrini \cite{Ajanalysemath1990, Aapplicable1991}. Indeed,  when $D = \emptyset$ and $\sigma$ is constant on $\mathbb R^N$,  he  considered the problem where the initial data in \eqref{heat initial} is replaced by the general data $u_0$ in problem \eqref{heat equation initial-boundary}-\eqref{heat initial}.  Then he  proved that if all the spatial isothermic surfaces of $u$ are stationary, then either $u_0-1$ is an eigenfunction of the Laplacian or $\Omega$ is a ball where $u_0$ is radially symmetric. See also \cite{Sjanalysemath1999, MMnonlinearAnal2016} for this direction.

The following sections are organized as follows. In section \ref{section2}, we give four preliminaries where the balance laws given in \cite{MSmathz1999, MSannals2002} play a key role on behalf of Varadhan's formula (see \eqref{varadhan formula}) given in \cite{Vcpam1967}. 
Section \ref{section3} is devoted to the proof of Theorem \ref{th:stationary isothermic}. Auxiliary functions $U, V$ given in section \ref{section3}  play a key role. If $D$ is not a ball, we use the transmission condition \eqref{transmission condition} on $\partial D$ to get a contradiction to Hopf's boundary point lemma. In section \ref{section4}, we prove Theorem \ref{th:stationary isothermic cauchy} by following the proof of Theorem \ref{th:stationary isothermic}. Auxiliary functions $U, V, W$  given in section \ref{section4}  play a key role. We notice that almost the same arguments work as in the proof of Theorem \ref{th:stationary isothermic}.

\setcounter{equation}{0}
\setcounter{theorem}{0}

\section{Preliminaries for $N \ge 2$}
\label{section2}

Concerning the behavior of the solutions of problem \eqref{heat equation initial-boundary}-\eqref{heat initial} and problem \eqref{heat Cauchy}, we start with the following lemma.
\begin{lemma} 
\label{le:initial behavior exponential decay and decay at infinity} Let $u$ be the solution of either problem \eqref{heat equation initial-boundary}-\eqref{heat initial} or problem \eqref{heat Cauchy}. We have the following assertions:
\begin{itemize}
\item[\rm (1)] For every compact set $K \subset \Omega$, there exist two positive constants $B$ and $b$ satisfying
$$
0 < u(x,t) < B e^{-\frac bt}\quad \mbox{ for every } (x,t) \in K \times (0,1].
$$
\item[\rm (2)] There exists a constant $M > 0$ satisfying
$$
0 \le 1-u(x,t) \le \min\{ 1, M t^{-\frac N2} |\Omega|\} \quad \mbox{ for every }(x,t) \in \Omega\times (0,+\infty) \mbox{ or } \in \mathbb R^N \times (0,\infty),
$$
 where $|\Omega|$ denotes the Lebesgue measure of the set $\Omega$.
 \item[\rm(3)] For the solution $u$ of problem \eqref{heat equation initial-boundary}-\eqref{heat initial}, there exist two positive constants $C$ and $\lambda$ satisfying
 $$
 0 \le 1-u(x,t) \le Ce^{-\lambda t}  \quad \mbox{ for every }(x,t) \in \Omega\times (0,+\infty).
 $$
 \item[\rm(4)] For the solution $u$ of problem \eqref{heat Cauchy} where $N \ge 3$, there exist two positive constants $\beta$ and $L$ satisfying
 $$
 \beta^{-1} |x|^{2-N} \le \int_0^\infty (1-u(x,t))\ dt \le \beta |x|^{2-N}\  \mbox{ if }\ |x| \ge L,
 $$
 where $\overline{\Omega} \subset B_L(0) = \{ x \in \mathbb R^N : |x| < L \}.$
 \end{itemize}
\end{lemma}

\noindent
{\it Proof.\ } We make use of the Gaussian bounds for the fundamental solutions of parabolic equations due to
Aronson\cite[Theorem 1, p. 891]{Ar1967bams}(see also \cite[p. 328]{FaS1986arma}). Let $g = g(x,t;\xi,\tau)$ be the fundamental solution of $u_t=\mbox{ div}(\sigma\nabla u)$. Then there exist two positive constants $\alpha$ and $M$ such that
\begin{equation}
\label{Gaussian bounds}
M^{-1}(t-\tau)^{-\frac N2}e^{-\frac{\alpha|x-\xi|^2}{t-\tau}}\le g(x,t;\xi,\tau) \le M(t-\tau)^{-\frac N2}e^{-\frac{|x-\xi|^2}{\alpha(t-\tau)}} 
\end{equation}
 for all $(x,t), (\xi,\tau) \in \mathbb R^N\times(0,+\infty)$ with $t > \tau$. 
 
 For the solution $u$ of problem \eqref{heat Cauchy},  $1-u$ is regarded as the unique bounded solution of the Cauchy problem for the diffusion equation with initial data ${\mathcal X}_{\Omega}$ which is greater
than or equal to the corresponding solution of the initial-boundary value problem for the diffusion equation under the homogeneous Dirichlet boundary condition by the comparison principle. Hence we have from \eqref{Gaussian bounds}
 $$
 1-u(x,t) = \int_{\mathbb R^N}  g(x,t;\xi,0){\mathcal X}_{\Omega}(\xi)\ d\xi \le M t^{-\frac N2}|\Omega|.
 $$
 The inequalities $0 \le 1-u \le 1$ follow  from the comparison principle. This completes the proof of (2). Moreover, (4) follows from \eqref{Gaussian bounds} as is noted in \cite[5. Remark, pp. 895--896]{Ar1967bams}.

For (1),  let $K$ be a compact set contained in $\Omega$.
 We set
$$
 \mathcal N_\rho = \{ x \in \mathbb R^N : \mbox{ dist}(x, \partial\Omega) < \rho \}
$$
 where $\rho = \frac 12 \mbox{ dist}(K, \partial\Omega)\ ( > 0)$. Define $v = v(x,t)$ by 
$$
 v(x,t) = \lambda \int_{\mathcal N_\rho} g(x,t; \xi,0)\ d\xi\quad \mbox{ for every }(x,t) \in \mathbb R^N\times(0,+\infty),
$$
 where a number $\lambda > 0$ will be determined later. Then it follows from \eqref{Gaussian bounds} that
$$
  v(x,t) \ge \lambda M^{-1} t^{-\frac N2}\int_{\mathcal N_\rho}e^{-\frac{\alpha|x-\xi|^2}{t}} d\xi \quad \mbox{ for } (x,t) \in \mathbb R^N\times(0,+\infty)
$$
 and hence we can choose $\lambda > 0$ satisfying
$$
  v \ge 1\ \mbox{ on } \partial\Omega \times (0,1].
$$
  Thus the comparison principle yields that
\begin{equation}
\label{upper bound by v}
  u \le v\ \mbox{ in } \Omega\times (0,1].
\end{equation}
On the other hand, it follows from \eqref{Gaussian bounds} that
$$
v(x,t) \le \lambda M t^{-\frac N2}\int_{\mathcal N_\rho}e^{-\frac{|x-\xi|^2}{\alpha t}} d\xi\quad \mbox{ for } (x,t) \in \mathbb R^N\times(0,+\infty).
$$
Since $|x-\xi| \ge \rho$ for every $x \in K$ and $\xi \in \mathcal N_\rho$, we observe that
$$
 v(x,t) \le \lambda M t^{-\frac N2}e^{-\frac {\rho^2}{\alpha t}}\left|\mathcal N_\rho\right| \quad \mbox{ for every }(x,t) \in K \times (0,+\infty),
 $$
 where $\left|\mathcal N_\rho\right|$ denotes the Lebesgue measure of the set $\mathcal N_\rho$.
 Therefore \eqref{upper bound by v} gives (1). 
 
  For (3),  for instance choose a large ball $B$ with $\overline{\Omega} \subset B$ and let $\varphi= \varphi(x)$ be the first positive eigenfunction of the  problem
  $$
  -\mbox{ div}(\sigma\nabla \varphi ) = \lambda \varphi \ \mbox{ in } B \ \mbox{ and } \varphi = 0\ \mbox{ on } \partial B
  $$
  with $\sup\limits_{B} \varphi= 1$. Choose $C > 0$ sufficiently large to have
  $$
  1 \le C \varphi\ \mbox{ in } \overline{\Omega}.
  $$
  Then it follows from the comparison principle that
  $$
  1-u(x,t) \le Ce^{-\lambda t} \varphi(x)\quad \mbox{ for every } (x,t) \in \Omega\times (0,+\infty),
  $$
  which gives (3). \qed
 
 The following asymptotic formula of the heat content of a ball touching at $\partial\Omega$ at only one point tells us about the interaction between the initial behavior of solutions and geometry of domain.
\begin{proposition} 
\label{prop:heat content asymptotics} Let $u$ be the solution of either problem \eqref{heat equation initial-boundary}-\eqref{heat initial} or problem \eqref{heat Cauchy}. 
Let $x \in \Omega$ and assume that the open ball $B_r(x)$  with radius $r >0$ centered at $x$
is contained in  $\Omega$ and such that $\overline{B_r(x)} \cap \partial\Omega = \{ y \}$ for some $y \in \partial\Omega$.
Then we have:
\begin{equation}
\label{asymptotics and curvatures}
\lim_{t\to +0}t^{-\frac{N+1}4 }\!\!\!\int\limits_{B_r(x)}\! u(z,t)\ dz=
C(N, \sigma)\left\{\prod\limits_{j=1}^{N-1}\left(\frac 1r - \kappa_j(y)\right)\right\}^{-\frac 12}.
\end{equation}
Here, 
$\kappa_1(y),\dots,\kappa_{N-1}(y)$ denote the principal curvatures of $\partial\Omega$ at $y$ with 
respect to the inward normal direction to $\partial\Omega$  
and $C(N, \sigma)$ is a positive constant given by
$$
C(N,\sigma) = \left\{\begin{array}{rll}2\sigma_s^{\frac {N+1}4}c(N) \ &\mbox{ for problem \eqref{heat equation initial-boundary}-\eqref{heat initial} },
\\
\frac {2\sqrt{\sigma_m}}{\sqrt{\sigma_s}+\sqrt{\sigma_m}}\sigma_s^{\frac {N+1}4}c(N) &\mbox{ for problem \eqref{heat Cauchy} },
\end{array}\right.
$$
where $c(N)$ is a positive constant depending only on $N$. (Notice that if $\sigma_s=\sigma_m$ then $C(N, \sigma) = \sigma_s^{\frac {N+1}4}c(N)$
for problem \eqref{heat Cauchy}, that is, just half of the constant  for problem \eqref{heat equation initial-boundary}-\eqref{heat initial}.)

\par
When $\kappa_j(y) = 1/r$ for some $j \in \{ 1, \cdots, N-1\}$, 
\eqref{asymptotics and curvatures} holds by setting the right-hand side to $+\infty$ (notice that 
$\kappa_j(y) \le 1/r$ always holds
for all $j$'s).
\end{proposition}

 \noindent
{\it Proof.\ } For the one-phase problem, that is, for the heat equation $u_t = \Delta u$, this lemma has been proved in
\cite[Theorem 1.1, p. 238]{MSjde2012} or in \cite[Theorem B, pp. 2024-2025 and Appendix, pp. 2029--2032]{MSmmas2013}.
The proof in \cite{MSmmas2013} was carried out by constructing appropriate super- and subsolutions in a neighborhood of $\partial\Omega$ in a short time with the aid of the initial behavior \cite[Lemma B.2, p. 2030]{MSmmas2013} obtained by Varadhan's formula \cite{Vcpam1967} for the heat equation $u_t = \Delta u$
\begin{equation}
\label{varadhan formula}
-4t \log u(x,t) \to \mbox{ dist}(x,\partial\Omega)^2\ \mbox{ as } t \to +0 \mbox{ uniformly on every compact set in $\Omega$.}
\end{equation}
 (See also \cite[Theorem A, p. 2024]{MSmmas2013} for the formula.)
Here, with no need of Varadhan's formula, (1) of Lemma \ref{le:initial behavior exponential decay and decay at infinity} gives sufficient information on 
 the initial behavior \cite[Lemma B.2, p. 2030]{MSmmas2013}.  We remark that since problem \eqref{heat equation initial-boundary}-\eqref{heat initial} is one-phase with conductivity $\sigma_s$ near $\partial\Omega$, we can obtain formula \eqref{asymptotics and curvatures} for problem \eqref{heat equation initial-boundary}-\eqref{heat initial} only by scaling in $t$. On the other hand, problem \eqref{heat Cauchy} is two-phase with conductivities $\sigma_m, \sigma_s$ near $\partial\Omega$
if $\sigma_m\not=\sigma_s$. Therefore, it is enough for us to prove formula \eqref{asymptotics and curvatures} for problem \eqref{heat Cauchy} where $\sigma_m\not=\sigma_s$. 

Let $u$ be the solution of problem \eqref{heat Cauchy} where $\sigma_m\not=\sigma_s$, and let us prove this lemma by modifying the proof of Theorem B in \cite[Appendix, pp. 2029--2032]{MSmmas2013}.

Let us consider the signed distance function $d^*= d^*(x)$ 
of $x \in \mathbb R^N$ to the boundary $\partial\Omega$ defined by
\begin{equation}
\label{signed distance}
d^*(x) = \left\{\begin{array}{rll}
 \mbox{ dist}(x,\partial\Omega)\ &\mbox{ if }\ x \in \Omega,
\\
-\mbox{ dist}(x,\partial\Omega)\ &\mbox{ if  }\ x \not\in \Omega.
\end{array}\right.
\end{equation}
Since $\partial\Omega$  is bounded and of class $C^2$,  there exists a number $\rho_0 > 0$ such that $d^*(x)$ is $C^2$-smooth on a compact neighborhood  $\mathcal N$ of the boundary $\partial\Omega$ given by
\begin{equation}
\label{neighborhood of boundary from both sides}
\mathcal N = \{ x \in \mathbb R^N : -\rho_0 \le d^*(x) \le \rho_0 \}.
\end{equation}
We make $\mathcal N$ satisfy $\mathcal N \cap \overline{D} = \emptyset$.
Introduce a function $F = F(\xi)$ for $\xi \in \mathbb R$ by
$$
F(\xi) = \frac 1{2\sqrt{\pi}} \int_\xi^\infty e^{-s^2/4} ds.
$$
Then $F$ satisfies
\begin{eqnarray*}
&&F^{\prime\prime} + \frac 12\xi F^\prime = 0\  \mbox{ and }  F^\prime < 0\ \mbox{ in } \mathbb R,
\\
&& F(-\infty) = 1,\ F(0) = \frac 12, \ \mbox{ and } F(+\infty) = 0.
\end{eqnarray*}
For each $\varepsilon \in (0,1/4)$, we define two functions $F_\pm = F_\pm(\xi)$ for $\xi \in \mathbb R$ by
$$
F_\pm(\xi) = F(\xi \mp 2\varepsilon).
$$
Then $F_\pm$ satisfies
\begin{eqnarray*}
&&F_\pm^{\prime\prime} + \frac 12\xi F_\pm^\prime = \pm\varepsilon F_\pm^\prime,\  F_\pm^\prime < 0\  \mbox{ and } F_- < F < F_+\ \mbox{ in } \mathbb R,
\\
&& F_\pm(-\infty) = 1,\ F_\pm(0) \gtrless   \frac 12, \ \mbox{ and } F_\pm(+\infty) = 0.
\end{eqnarray*}

By setting $\eta = t^{-\frac 12}d^*(x),\ \mu =  {\sqrt{\sigma_m}}/{\sqrt{\sigma_s}}$ and $\theta_\pm = 1 + (\mu-1)F_\pm(0)\ ( > 0)$, we introduce two functions $v_{\pm} = v_{\pm}(x,t)$ by
\begin{equation}
\label{def-of-pre-subsupersolutions}
v_{\pm}(x,t) =\left\{\begin{array}{rl} \frac {\mu}{\theta_\pm}F_{\pm}\left(\sigma_s^{-\frac 12}\eta\right) &\mbox{ for }\ (x,t) \in  \Omega \times (0,+\infty),
\\
\frac {1}{\theta_\pm}\left[F_{\pm}\left(\sigma_m^{-\frac 12}\eta\right) + \theta_\pm-1\right] &\mbox{ for }\ (x,t) \in \left(\mathbb R^N\setminus \Omega\right) \times (0,+\infty).
\end{array}\right.
\end{equation}
Then $v_{\pm}$ satisfies the transmission conditions 
\begin{equation}
\label{transmission conditions for pre-subsupersolutions}
v_{\pm}\big|_+ = v_{\pm}\big|_-\ \mbox{ and }\ \sigma_m\frac {\partial v_{\pm}}{\partial\nu}\Big|_+ = \sigma_s\frac {\partial v_{\pm}}{\partial\nu}\Big|_-\ \mbox{ on } \partial\Omega\times (0,+\infty),
\end{equation}
where $+$ denotes the limit from outside and $-$ that from inside of $\Omega$ and $\nu=\nu(x)$ denotes the outward unit normal vector to $\partial\Omega$ at $x \in \partial\Omega$, since $\nu = -\nabla d^*$ on $\partial\Omega$. Moreover we observe that
for each $\varepsilon \in (0,1/4),$ there exists $t_{1,\varepsilon} \in (0,1]$ satisfying
 \begin{equation}
 \label{differential inequalities for vpm}
 (\pm1)\left\{(v_\pm)_t - \sigma\Delta v_\pm\right\} > 0\quad\mbox{ in }\ \left(\mathcal N\setminus\partial\Omega\right) \times (0,t_{1, \varepsilon}].
 \end{equation}
In fact, a straightforward computation gives   
$$
 (v_\pm)_t - \sigma\Delta v_\pm =\left\{\begin{array}{rl}
  -\frac{\mu}{t\theta_\pm} \left(\pm\varepsilon +\sqrt{\sigma_s t}\Delta d^*\right)\ F_\pm^\prime \quad&\mbox{ in }\ \left(\mathcal N\cap\Omega\right) \times (0,+\infty),
  \\
 -\frac{1}{t\theta_\pm} \left(\pm\varepsilon +\sqrt{\sigma_m t}\Delta d^*\right)\ F_\pm^\prime \quad&\mbox{ in }\ \left(\mathcal N\setminus\overline{\Omega}\right) \times (0,+\infty).
\end{array}\right. 
 $$
Then, for each $\varepsilon \in (0,1/4)$, by setting
$
t_{1,\varepsilon} = \frac 1{\max\{\sigma_s,\sigma_m\}} \left(\frac \varepsilon{2M}\right)^2,
$
where $M = \max\limits_{x \in \mathcal N} |\Delta d^*(x)|$, we obtain \eqref{differential inequalities for vpm}.

Then, in view of (1) of Lemma \ref{le:initial behavior exponential decay and decay at infinity} and the definition (\ref{def-of-pre-subsupersolutions}) of $v_\pm,$ we see that
there exist two positive constants $E_1$ and $E_2$ satisfying
\begin{equation}
\label{lm:estimate-inside}
\max\{ |v_+|,\ |v_-|,\ |u| \} \le E_1 e^{-\frac {E_2}t}\ \mbox{ in }\ \overline{\Omega\setminus\mathcal N} \times (0,1].
\end{equation}

\par
By setting, for $(x,t) \in \mathbb R^N \times (0,+\infty)$,
\begin{equation}
\label{super- and subsolutions}
w_\pm(x,t)=(1\pm\varepsilon) v_{\pm}(x,t)\pm 2E_1 e^{-\frac {E_2}t},
\end{equation}
since $v_\pm$ and $u$ are all nonnegative, we obtain from \eqref{lm:estimate-inside} that
\begin{equation}
\label{inequality inside compact set}
 w_-\le u \le w_+\ \mbox{ in } \overline{\Omega\setminus\mathcal N} \times (0,1].
 \end{equation}
Moreover, in view of the facts that $F_\pm(-\infty) = 1$ and $F_\pm(+\infty) = 0$, we see  that
there exists $t_\varepsilon \in (0, t_{1,\varepsilon}]$ satisfying
\begin{equation}
\label{initial and boundary conditions are OK}
w_-\le u \le w_+\ \mbox{ on } \left( (\partial\mathcal N\setminus\Omega) \times (0,t_\varepsilon] \right) \cup \left(\mathcal N \times \{0\}\right).
\end{equation}
Then, in view of \eqref{transmission conditions for pre-subsupersolutions}, \eqref{differential inequalities for vpm},  \eqref{inequality inside compact set}, \eqref{initial and boundary conditions are OK} and the definition \eqref{super- and subsolutions} of $w_\pm$,  
we have from the comparison principle over $\mathcal N$  that
\begin{equation}
\label{estimates from both}
w_- \le u \le w_+\quad\mbox{ in }\ \left(\overline{\mathcal N} \cup \Omega\right) \times(0,t_\varepsilon].
\end{equation}

By writing
$$
\Gamma_s = \{ x \in \Omega : d^*(x) = s \}\ \mbox{ for } s > 0,
$$
let us quote a geometric lemma from \cite{MSedinburgh2007} adjusted to our situation.
\begin{lemma}{\rm (\cite[Lemma 2.1, p. 376]{MSedinburgh2007})}
\label{lm:asympvol} If $\displaystyle{ \max_{1\le j \le N-1} \kappa_j(y) < \frac 1r}$, then
we have:
\begin{equation*}
\label{asympvol}
\lim_{s\to 0^+} s^{-\frac{N-1}{2}} \mathcal H^{N-1}(\Gamma_s\cap B_r(x))=2^{\frac{N-1}{2}}\omega_{N-1} \
\left\{\prod_{j=1}^{N-1}\left(\frac1{r}-\kappa_j(y)\right)\right\}^{-\frac12},
\end{equation*}
where $\mathcal H^{N-1}$ is the standard $(N-1)$-dimensional Hausdorff measure, and $\omega_{N-1}$ is the volume of the unit ball in $\mathbb R^{N-1}.$
\end{lemma}

Let us consider the case where $\displaystyle{ \max_{1\le j \le N-1} \kappa_j(y) < \frac 1r}$.
 Then it follows from \eqref{estimates from both}  that for every $t \in (0,t_\varepsilon]$
\begin{equation}
\label{estimates from both sides}
t^{-\frac{N+1}4 }\int_{B_r(x)} w_- \ dz \le t^{-\frac{N+1}4 } \int_{B_r(x)} u \ dz \le t^{-\frac{N+1}4 }\int_{B_r(x)} w_+ \ dz.
\end{equation}
On the other hand, with the aid of the co-area formula, we have:
\begin{eqnarray*}
&&\int_{B_r(x)} v_\pm\ dz = 
\\
&&\frac \mu{\theta_\pm}(\sigma_st)^{\frac{N+1}4}\int_0^{2r(\sigma_st)^{-\frac 12}} F_\pm(\xi)\xi^{\frac {N-1}2} \left((\sigma_st)^{\frac 12}\xi\right)^{-\frac {N-1}2} \mathcal H^{N-1}\left(\Gamma_{(\sigma_st)^{\frac12}\xi}\cap B_r(x)\right) d\xi,
\end{eqnarray*}
where $v_\pm$ is defined by \eqref{def-of-pre-subsupersolutions}. Thus,
by Lebesgue's dominated convergence theorem and Lemma \ref{lm:asympvol},  we get
\begin{eqnarray*}
&&\lim_{t\to +0} t^{-\frac{N+1}4}\int_{B_r(x)} w_\pm\ dx =
\\
&& \frac \mu{\theta_\pm}(\sigma_s)^{\frac{N+1}4}2^{\frac{N-1}2}\omega_{N-1}\left\{\prod_{j=1}^{N-1}\left(\frac 1r-\kappa_j(y)\right)\right\}^{-\frac 12}\int_0^\infty  F_\pm(\xi) \xi^{\frac{N-1}2} d\xi.
\end{eqnarray*}
Moreover, again by Lebesgue's dominated convergence theorem, since 
$$
\lim\limits_{\varepsilon \to 0}\theta_\pm = 1 + (\mu-1)F(0) = \frac {\mu+1}2\ \mbox{ and } \mu =  {\sqrt{\sigma_m}}/{\sqrt{\sigma_s}},
$$
we see  that
\begin{eqnarray*}
&&\lim_{t\to +0} t^{-\frac{N+1}4}\int_{B_r(x)} w_\pm\ dx =
\\
&& \frac {2\sqrt{\sigma_m}}{\sqrt{\sigma_s}+\sqrt{\sigma_m}}(\sigma_s)^{\frac{N+1}4}2^{\frac{N-1}2}\omega_{N-1}\left\{\prod_{j=1}^{N-1}\left(\frac 1r-\kappa_j(y)\right)\right\}^{-\frac 12}\int_0^\infty  F(\xi) \xi^{\frac{N-1}2} d\xi.
\end{eqnarray*}
Therefore \eqref{estimates from both sides} gives formula \eqref{asymptotics and curvatures} provided
$\displaystyle{ \max_{1\le j \le N-1} \kappa_j(y) < \frac 1r}$.

Once this is proved, the case where $\kappa_j(y) = 1/r$ for some $j \in \{ 1, \cdots, N-1\}$ can be dealt with as in \cite[p. 248]{MSjde2012} by choosing a sequence of balls $\{ B_{r_k}(x_k) \}_{k=1}^\infty$  satisfying:
$$
 r_k < r,\ y \in \partial B_{r_k}(x_k), \mbox{ and } B_{r_k}(x_k) \subset B_r(x)\mbox{ for every } k \ge 1,\ \mbox{ and }\ \lim_{k\to \infty}r_k = r.
$$
Then, because of $\displaystyle{ \max_{1\le j \le N-1} \kappa_j(y) \le \frac 1r < \frac 1{r_k}}$, applying formula \eqref{asymptotics and curvatures} to each ball $B_{r_k}(x_k)$  yields that
$$
\liminf_{t \to +0} t^{-\frac{N+1}4}\int_{B_r(x)} u(z,t)\ dz = +\infty.
$$
This completes the proof of Proposition \ref{prop:heat content asymptotics}. \qed

In order to determine the symmetry of $\Omega$, we employ the following lemma.
\begin{lemma} 
\label{le: constant weingarten curvature}
Let $u$ be the solution of either problem \eqref{heat equation initial-boundary}-\eqref{heat initial} or problem \eqref{heat Cauchy}. 
Under the assumption \eqref{stationary isothermic surface partially} of {\rm Theorem \ref{th:stationary isothermic}} and {\rm Theorem \ref{th:stationary isothermic cauchy}}, the following assertions hold:
\begin{enumerate}[\rm (1)]
\item There exists a number $R > 0$ such that 
$$
\mbox{\rm dist}(x, \partial\Omega) = R\ \mbox{ for every } x \in \Gamma;
$$
\item $\Gamma$ is a real analytic hypersurface;
\item there exists a connected component $\gamma$ of $\partial\Omega$, that is also a real analytic hypersurface, such that the mapping $\gamma \ni y \mapsto x(y)  \equiv y-R\nu(y) \in \Gamma$, where $\nu(y)$ is the outward unit normal vector to $\partial\Omega$ at $y \in \gamma$,  is a diffeomorphism; in particular $\gamma$ and $\Gamma$  are parallel hypersurfaces at distance $R$;
\item it holds that
\begin{equation}
\label{bounds of curvatures}
 \max_{1\le j \le N-1}\kappa_j(y) < \frac 1R\ \mbox{ for every } y \in \gamma,
\end{equation}
where $\kappa_1(y), \cdots, \kappa_{N-1}(y)$ are the principal curvatures of $\partial\Omega$ at $y \in \gamma$ with respect to the inward unit normal vector $-\nu(y)$ to $\partial\Omega$;
\item there exists a number $c > 0$ such that
\begin{equation}
\label{monge-ampere}
\prod_{j=1}^{N-1} \left(\frac 1R-\kappa_j(y)\right) = c\quad\mbox{ for every } y \in \gamma.
\end{equation}
\end{enumerate}
\end{lemma}
 
\noindent
{\it Proof.\ } First it follows from  the assumption \eqref{near the boundary} that 
$$
B_r(x) \subset \Omega\setminus\overline{D}\ \mbox{ for every } x \in \partial G \mbox{ with } 0 < r \le \mbox{ dist}(x,\partial\Omega).
$$
Therefore, since $\sigma = \sigma_s$ in $\Omega\setminus\overline{D}$, we can use a balance law (see \cite[Theorem 2.1, pp. 934--935]{MSannals2002} or \cite[Theorem 4, p. 704]{MSmathz1999}) to obtain from \eqref{stationary isothermic surface partially} that
for every $p, q \in \Gamma$ and $t > 0$
\begin{equation}
\label{balance law}
\int_{B_r(p)}u(z,t)\ dz = \int_{B_r(q)}u(z,t)\ dz\ \mbox{ if }\  0 < r \le \min\{\mbox{ dist}(p,\partial\Omega),\mbox{ dist}(q,\partial\Omega) \}.  
\end{equation}
Let us show assertion (1). Suppose that there exist a pair of points $p$ and $q$ satisfying
$$
\mbox{ dist}(p,\partial\Omega) < \mbox{ dist}(q,\partial\Omega).
$$
Set $r = \mbox{ dist}(p,\partial\Omega)$. Then there exists a point $y \in \partial\Omega$ such that
$y \in \overline{B_r(p)}\cap\partial\Omega$. Choose a smaller ball $B_{\hat{r}}(x) \subset B_r(p)$ with $0 < \hat{r} < r$ and 
$\overline{B_{\hat{r}}(x)} \cap \partial B_r(p) = \{ y \}$. Since $\displaystyle{ \max_{1\le j \le N-1} \kappa_j(y) \le \frac 1r < \frac 1{\hat{r}}}$, by applying Proposition \ref{prop:heat content asymptotics} to the ball $B_{\hat{r}}(x)$, we get
$$
\liminf_{t \to +0}t^{-\frac{N+1}4 }\int_{B_r(p)} u(z,t)\ dz \ge \lim_{t \to +0}t^{-\frac{N+1}4 }\int_{B_{\hat{r}}(x)} u(z,t)\ dz > 0.
$$
On the other hand, since $\overline{B_r(q)} \subset \Omega$, it follows from (1) of Lemma \ref{le:initial behavior exponential decay and decay at infinity} that
$$
\lim_{t \to +0}t^{-\frac{N+1}4 }\int_{B_r(q)} u(z,t)\ dz = 0,
$$
which contradicts \eqref{balance law}, and hence assertion (1) holds true.

We can find a point $x_* \in \Gamma$ and a ball $B_\rho(z_*)$ such that $B_\rho(z_*) \subset G$ and $x_* \in \partial B_\rho(z_*)$. Since $\Gamma$ satisfies \eqref{nearest component}, assertion (1) yields that there exists a point $y_* \in \partial\Omega$ satisfying
$$
B_{R+\rho}(z_*) \subset \Omega,\  y_* \in \overline{B_{R+\rho}(z_*)} \cap \partial\Omega, \mbox{ and } \overline{B_R(x_*)} \cap \partial\Omega = \{y_*\}.
$$
Observe that 
$$
\displaystyle{ \max_{1\le j \le N-1} \kappa_j(y_*)\le \frac 1{R+\rho} < \frac 1R}\ \mbox{ and }\ x_* = y_*-R\nu(y_*) \equiv x(y_*).
$$
 Define $\gamma \subset \partial\Omega$ by
$$
\gamma = \left\{ y \in \partial\Omega : \overline{B_R(x)}\cap\partial\Omega =\{y\} \ \mbox{ for } x = y-R\nu(y) \in \Gamma \mbox{ and }
\max_{1\le j \le N-1} \kappa_j(y) < \frac 1R \right\}.
$$
Hence $y_* \in \gamma$ and $\gamma \not=\emptyset$.  By Proposition \ref{prop:heat content asymptotics} we have that for every $y \in \gamma$ and $x = x(y) (=y-R\nu(y))$
\begin{equation}
\label{key asymptotics for the proof}
\lim_{t \to +0}t^{-\frac{N+1}4 }\int_{B_R(x)} u(z,t)\ dz = C(N,\sigma)\left\{\prod\limits_{j=1}^{N-1}\left(\frac 1R - \kappa_j(y)\right)\right\}^{-\frac 12}.
\end{equation}
Here let us show that, if $y \in \gamma$ and $x = x(y)$, then $\nabla u(x,t) \not= 0$ for some $t > 0$, which guarantees that
in a neighborhood of $x$, $\Gamma$ is a part of a real analytic hypersurface properly embedded in $\mathbb R^N$ because of \eqref{stationary isothermic surface partially},
real analyticity of $u$ with respect to the space variables, and the implicit function theorem. Moreover, this together with the implicit function theorem guarantees that $\gamma$ is open in $\partial\Omega$ and the mapping 
$\gamma\ni y\mapsto x(y) \in \Gamma$ is a local diffeomorphism, which is also real analytic. If we can prove additionally that $\gamma$ is closed in $\partial\Omega$, then the mapping $\gamma\ni y\mapsto x(y) \in \Gamma$ is a diffeomorphism and $\gamma$ is a connected component of $\partial\Omega$ since $\Gamma$ is a connected component of $\partial G$, and hence all the remaining assertions (2) -- (5) follow from \eqref{balance law}, \eqref{key asymptotics for the proof} and the definition of $\gamma$. We shall prove this later in the end of the proof of Lemma \ref{le: constant weingarten curvature}. 

Before this we show that, if $y \in \gamma$ and $x = x(y)$, then $\nabla u(x,t) \not= 0$ for some $t > 0$.
Suppose that $\nabla u(x,t) = 0$ for every $t > 0$. Then we use another balance law (see \cite[Corollary 2.2, pp. 935--936]{MSannals2002}) to obtain that
\begin{equation}
\label{balance for the gradient}
\int_{B_R(x)} (z-x) u(z,t)\ dz = 0\ \mbox{ for every } t > 0.
\end{equation}
On the other hand, (1) of Lemma \ref{le:initial behavior exponential decay and decay at infinity} yields that
\begin{equation}
\label{compact sets are negligible initially}
\lim_{t \to +0}t^{-\frac{N+1}4 }\int_K u(z,t)\ dz =0\ \mbox{ for every compact set }K \subset \Omega,
\end{equation}
and hence by \eqref{key asymptotics for the proof} it follows that for every $\varepsilon > 0$
\begin{equation}
\label{small neighborhood of tangent point is enough}
\lim_{t \to +0}t^{-\frac{N+1}4 }\int_{B_R(x)\cap B_\varepsilon(y)} u(z,t)\ dz = C(N,\sigma)\left\{\prod\limits_{j=1}^{N-1}\left(\frac 1R - \kappa_j(y)\right)\right\}^{-\frac 12}.
\end{equation}
This implies that
$$
\lim_{t \to +0}t^{-\frac{N+1}4 }\int_{B_R(x)} (z-x) u(z,t)\ dz = C(N,\sigma)\left\{\prod\limits_{j=1}^{N-1}\left(\frac 1R - \kappa_j(y)\right)\right\}^{-\frac 12}(y-x) \not=0,
$$
which contradicts \eqref{balance for the gradient}.

It remains to show that $\gamma$ is closed in $\partial\Omega$. Let $\{ y^n \}$ be a sequence of points in $\gamma$ with 
$\lim\limits_{n \to \infty} y^n = y^\infty \in \partial\Omega$, and let us prove that $y^\infty \in \gamma$. By combining \eqref{balance law} with \eqref{key asymptotics for the proof}, we see that there exists a positive number $c$ satisfying assertion (5) and hence by continuity
\begin{equation}
\label{monge ampere equation we wanted}
\prod_{j=1}^{N-1} \left(\frac 1R-\kappa_j(y^\infty)\right) = c > 0\ \mbox{ and } \max_{1\le j \le N-1} \kappa_j(y^\infty) \le \frac 1R, 
\end{equation}
since $y^j \in \gamma$ for every $j$. Thus $\displaystyle{\max_{1\le j \le N-1} \kappa_j(y^\infty) < \frac 1R}$. Let 
$x^\infty = y^\infty -R\nu(y^\infty) (=x(y^\infty)).$ It suffices to show that $\overline{B_R(x^\infty)}\cap\partial\Omega =\{y^\infty\}$. Suppose that there exists another point $y \in \overline{B_R(x^\infty)}\cap\partial\Omega$. Then for every $\hat{R} \in (0, R)$ we can find two points $p^\infty$ and $p$ in $B_R(x^\infty)$ such that
$$
B_{\hat{R}}(p^\infty) \cup B_{\hat{R}}(p) \subset B_R(x^\infty), \ \overline{B_{\hat{R}}(p^\infty)} \cap \partial\Omega = \{y^\infty\},\ \mbox{ and } \overline{B_{\hat{R}}(p)} \cap \partial\Omega = \{y\}.
$$
Hence by Proposition \ref{prop:heat content asymptotics} we have
\begin{eqnarray*}
&&\lim_{t \to +0}t^{-\frac{N+1}4 }\int_{B_{\hat{R}}(p^\infty)} u(z,t)\ dz = C(N,\sigma)\left\{\prod\limits_{j=1}^{N-1}\left(\frac 1{\hat{R}} - \kappa_j(y^\infty)\right)\right\}^{-\frac 12},
\\
&&\lim_{t \to +0}t^{-\frac{N+1}4 }\int_{B_{\hat{R}}(p)} u(z,t)\ dz = C(N,\sigma)\left\{\prod\limits_{j=1}^{N-1}\left(\frac 1{\hat{R}} - \kappa_j(y)\right)\right\}^{-\frac 12}.
\end{eqnarray*}
Thus, with the same reasoning as in \eqref{small neighborhood of tangent point is enough} by choosing small $\varepsilon > 0$, we have from \eqref{monge ampere equation we wanted}, \eqref{balance law}, \eqref{key asymptotics for the proof} and assertion (5) that  for every $x \in \gamma$
\begin{eqnarray*}
&& C(N,\sigma)\left\{\prod_{j=1}^{N-1} \left(\frac 1R-\kappa_j(y^\infty)\right)\right\}^{-\frac 12} = C(N,\sigma)c^{-\frac 12} 
\\
&&= \lim_{t \to +0}t^{-\frac{N+1}4 }\int_{B_{R}(x)} u(z,t)\ dz = \lim_{t \to +0}t^{-\frac{N+1}4 }\int_{B_{R}(x^\infty)} u(z,t)\ dz 
\\
&& \ge \lim_{t \to +0}t^{-\frac{N+1}4 }\int_{B_{\hat{R}}(p^\infty)\cap B_\varepsilon(y^\infty)} u(z,t)\ dz + \lim_{t \to +0}t^{-\frac{N+1}4 }\int_{B_{\hat{R}}(p)\cap B_\varepsilon(y)} u(z,t)\ dz
\\
&& =C(N,\sigma)\left[\left\{\prod_{j=1}^{N-1} \left(\frac 1{\hat{R}}-\kappa_j(y^\infty)\right)\right\}^{-\frac 12} + \left\{\prod_{j=1}^{N-1} \left(\frac 1{\hat{R}}-\kappa_j(y)\right)\right\}^{-\frac 12}\right].
\end{eqnarray*}
Since $\hat{R} \in (0, R)$ is arbitrarily chosen, this gives a contradiction, and hence $\gamma$ is closed in $\partial\Omega$. \qed

\begin{lemma} 
\label{le: constant weingarten curvature 2}
Let $u$ be the solution of  problem \eqref{heat Cauchy}. 
Under the assumption \eqref{stationary isothermic surface} of {\rm Theorem \ref{th:stationary isothermic cauchy}}, the same assertions {\rm (1)--(5)} as in {\rm Lemma \ref{le: constant weingarten curvature}} hold provided $\Gamma$ and $\gamma$ are replaced by $\partial G$ and $\partial\Omega$, respectively.
\end{lemma}

\noindent
{\it Proof.\ } 
By the same reasoning as in assertion (1) of Lemma \ref{le: constant weingarten curvature} we have assertion (1) from the assumption \eqref{stationary isothermic surface}. Since every component $\Gamma$ of $\partial G$ has the same distance $R$ to $\partial\Omega$, every component $\Gamma$ satisfies the assumption \eqref{nearest component}. Therefore, we can use the same arguments as in the proof of Lemma \ref{le: constant weingarten curvature} to prove this lemma. Here we must have
$$
\partial\Omega =\{ x \in \mathbb R^N : \mbox{ dist}(x, \overline{G}) = R \}. \ \mbox{ \qed }
$$

\setcounter{equation}{0}
\setcounter{theorem}{0}

\section{Proof of Theorem \ref{th:stationary isothermic}}
\label{section3}

Let $u$ be the solution of  problem \eqref{heat equation initial-boundary}-\eqref{heat initial} for $N \ge 2$.
With the aid of Aleksandrov's sphere theorem \cite[p. 412]{Alek1958vestnik},  Lemma \ref{le: constant weingarten curvature} yields that $\gamma$ and $\Gamma$ are concentric spheres. Denote by $x_0 \in \mathbb R^N$ the common center of $\gamma$ and $\Gamma$. By combining the initial and boundary conditions of problem \eqref{heat equation initial-boundary}-\eqref{heat initial} and  the assumption \eqref{stationary isothermic surface partially} with the real analyticity in $x$ of $u$ over $\Omega\setminus\overline{D}$,  we see that $u$ is radially symmetric with respect to $x_0$  in $x$ on $\left(\Omega\setminus\overline{D}\right) \times (0,\infty)$. Here we used the assumption that $\Omega\setminus\overline{D}$ is connected. Moreover, in view of the Dirichlet boundary condition \eqref{heat Dirichlet}, we can distinguish the following two cases:
$$
\mbox{\rm (I)  } \Omega \mbox{ is a ball;}\qquad \mbox{\rm (II)  } \Omega \mbox{ is a spherical shell.}
$$ 

By virtue of (3) of Lemma \ref{le:initial behavior exponential decay and decay at infinity},  we can introduce the following two auxiliary functions $U = U(x), \ V = V(x)$ by
\begin{eqnarray}
&& U(x) = \int_0^\infty(1-u(x,t) )\ dt\quad \mbox{ for } x \in \Omega \setminus\overline{D}, \label{auxiliary function on shell}
\\
&& V(x) = \int_0^\infty(1-u(x,t) )\ dt\quad \mbox{ for } x \in D.\label{auxiliary function on core}
\end{eqnarray}
Then we observe that
\begin{eqnarray}
&&-\Delta U = \frac 1{\sigma_s} \ \mbox{ in } \Omega \setminus\overline{D}, \ -\Delta V = \frac 1{\sigma_c} \ \mbox{ in } D,\label{poisson equations}
\\
&&U = V\ \mbox{ and }\  \sigma_s \frac {\partial U}{\partial\nu} = \sigma_c \frac {\partial V}{\partial\nu}  \ \mbox{ on } \partial D, \label{transmission condition}
\\
&& U = 0\ \mbox{ on } \partial\Omega, \label{Dirichlet conditions}
\end{eqnarray}
where $\nu = \nu(x)$ denotes the outward unit normal vector to $\partial D$ at $x \in \partial D$ and \eqref{transmission condition} is the transmission condition.
Since $U$ is radially symmetric with respect to $x_0$, by setting $r = |x-x_0|$ for $x \in  \Omega \setminus\overline{D}$ we have
\begin{equation}
\label{ODE radial}
-\frac {\partial^2}{\partial r^2}U - \frac {N-1}r \frac {\partial}{\partial r} U =  \frac 1{\sigma_s}\ \mbox{ in }  \Omega \setminus\overline{D}.
\end{equation}
Solving this ordinary differential equation yields that
\begin{equation}
U = \left\{\begin{array}{rll}
c_1 r^{2-N}- \frac 1{2N\sigma_s}r^2 + c_2 \ &\mbox{ if }\ N \ge 3,
\\
- c_1\log r  - \frac 1{4\sigma_s}r^2 + c_2 \ &\mbox{ if }\ N =2,
\end{array}\right.\label{ODE is solved}
\end{equation}
where $c_1, c_2$ are some constants depending on $N$. Remark that $U$ can be extended as a radially symmetric function of $r$ in $\mathbb R^N \setminus \{ x_0 \}$.

Let us first show that case (II) does not occur.  Set $\Omega = B_{\rho_+}(x_0) \setminus \overline{B_{\rho_-}(x_0)}$ for some numbers $\rho_+ > \rho_- > 0$.
Since $\Omega\setminus\overline{D}$ is connected, \eqref{Dirichlet conditions} yields that $U(\rho_+) = U(\rho_-) = 0$ and hence $c_1 < 0$. Moreover we observe that
\begin{equation}
U^{\prime\prime} < 0\ \mbox{ on } [\rho_- , \rho_+].\label{concave}
\end{equation}
Recall that $D$ may have finitely many connected components. Let us take a connected component $D_* \subset D$.
Then, since $\overline{D_*} \subset \Omega$, we see that there exist $\rho_* \in  (\rho_- , \rho_+)$  and $x_* \in \partial D_*$ which satisfy
\begin{equation}
\label{key point for Hopf 1}
U(\rho_*) = \min\{ U(r) : r = |x-x_0|, x \in \partial D_*\}\ \mbox{ and } \rho_* = |x_* - x_0|.
\end{equation}
Notice that $\nu(x_*)$ equals either $\frac{x_*-x_0}{\rho_*}$ or $- \frac{x_*-x_0}{\rho_*}$.
For $r > 0$, set
\begin{equation}
\label{key function 1}
\hat{U}(r) = U(\rho_*) + \frac{\sigma_s}{\sigma_c}(U(r)- U(\rho_*)).
\end{equation}
Since
\begin{equation}
\label{useful identity}
\hat{U}(r) - U(r) = \left( \frac {\sigma_s}{\sigma_c} - 1\right)(U(r) - U(\rho_*)),
\end{equation}
it follows that
\begin{equation}
\label{hat against the original}
\hat{U} \left\{\begin{array}{rll}
\ge U \ &\mbox{ if }\ \sigma_s > \sigma_c
\\
\le U  \ &\mbox{ if }\ \sigma_s < \sigma_c
\end{array}\right. \mbox{ on } \partial D_*.
\end{equation}
Moreover, we remark that $\hat{U}$ never equals $U$ identically on $\partial D_*$ since $\Omega\setminus\overline{D_*}$ is connected and $\Omega$ is a spherical shell. Observe that  
\begin{equation}
\label{hat function is good}
-\Delta \hat{U} = \frac 1{\sigma_c}\  \mbox{ and } \ \frac {\partial\hat{U}}{\partial r}= \frac {\sigma_s}{\sigma_c} \frac {\partial U}{\partial r}\ \mbox{ in } \overline{D_*}.
\end{equation}
On the other hand, we have
\begin{equation}
\label{the equations V satisfies}
-\Delta V = \frac 1{\sigma_c}\ \mbox{ in } D_*\  \mbox{ and } \ V = U\ \mbox{ on } \partial D_*.
\end{equation}
Then it follows from \eqref{hat against the original} and the strong comparison principle that
\begin{equation}
\label{strong comparison}
\hat{U} \left\{\begin{array}{rll}
> V \ &\mbox{ if }\ \sigma_s > \sigma_c
\\
< V \ &\mbox{ if }\ \sigma_s < \sigma_c
\end{array}\right. \mbox{ in } D_*,
\end{equation}
since $\hat{U}$ never equals $U$ identically on $\partial D_*$. The transmission condition \eqref{transmission condition} with the definition of $\hat{U}$ tells us that 
\begin{equation}
\label{transmission at the point}
\hat{U} = V\ \mbox{ and } \ \frac {\partial \hat{U}}{\partial\nu} = \frac {\partial V}{\partial\nu} \ \mbox{ at } x = x_* \in \partial D_*,
\end{equation}
since $\nu(x_*)$ equals either $\frac{x_*-x_0}{\rho_*}$ or $- \frac{x_*-x_0}{\rho_*}$.
Therefore applying  Hopf's boundary point lemma to the harmonic function $\hat{U} - V$ gives a contradiction to \eqref{transmission at the point},  and hence case (II) never occurs. (See \cite[Lemma 3.4, p. 34]{GT1983} for Hopf's boundary point lemma.)

Let us consider case (I).  Set $\Omega = B_{\rho}(x_0)$ for some number $\rho > 0$. We distinguish the following three cases:
$$
\mbox{\rm (i)  } c_1 = 0;\qquad \mbox{\rm (ii)  } c_1 > 0; \qquad \mbox{\rm (iii)  } c_1 < 0.
$$ 
We shall show that only case (i) occurs. Let us consider case (i) first. Note that
\begin{equation}
\label{monotone and smooth}
U^\prime(r) < 0 \ \mbox{ if } r > 0, \mbox{ and } U^\prime(0)= 0.
\end{equation}
Take an arbitrary component $D_* \subset D$. Then, since $\overline{D_*} \subset \Omega = B_{\rho}(x_0)$, we see that there exist $\rho_* \in  (0, \rho)$  and $x_* \in \partial D_*$ which also satisfy \eqref{key point for Hopf 1}. Notice that $\nu(x_*)$ equals $\frac{x_*-x_0}{\rho_*}$. For $r \ge 0$, define $\hat{U} = \hat{U}(r)$ by \eqref{key function 1}. Then, by \eqref{useful identity} we also have  \eqref{hat against the original}. Observe that both \eqref{hat function is good} and \eqref{the equations V satisfies} also hold true.
 Then it follows from \eqref{hat against the original} and the comparison principle that
\begin{equation}
\label{weak comparison}
\hat{U} \left\{\begin{array}{rll}
\ge V \ &\mbox{ if }\ \sigma_s > \sigma_c
\\
\le V \ &\mbox{ if }\ \sigma_s < \sigma_c
\end{array}\right. \mbox{ in } D_*.
\end{equation}
The transmission condition \eqref{transmission condition} with the definition of $\hat{U}$ also yields \eqref{transmission at the point} since $\nu(x_*)$ equals $\frac{x_*-x_0}{\rho_*}$. Therefore, by applying  Hopf's boundary point lemma to the harmonic function $\hat{U} - V$, we conclude from \eqref{transmission at the point} that
$$
\hat{U} \equiv V\ \mbox{ in } D_*
$$
and hence  $D_*$ must be a ball centered at $x_0$. In conclusion, $D$ itself is connected and must be a ball centered at $x_0$, 
since $D_*$ is an arbitrary component of $D$.

Next, let us show that case (ii) does not occur. In case (ii) we have
\begin{equation}
\label{monotone and +infinity}
U^\prime(r) < 0 \ \mbox{ if } r > 0, \ \lim_{r \to 0}U(r)= +\infty, \mbox{ and }\ x_0 \in D.
\end{equation}
Let us choose the connected component $D_*$ of $D$ satisfying $x_0 \in D_*$.  Then, since $\overline{D_*} \subset \Omega = B_{\rho}(x_0)$, we see that there exist $\rho_{*1},  \rho_{*2} \in  (0, \rho)$  and $x_{*1}, x_{*2} \in \partial D_*$ which  satisfy that $\rho_{*1} \le \rho_{*2}$ and
\begin{eqnarray}
&&U(\rho_{*1}) = \max\{ U(r) : r = |x-x_0|, x \in \partial D_*\} \mbox{ and } \rho_{*1} =|x_{*1}-x_0|, \label{maximum point}
\\
&&U(\rho_{*2}) = \min\{ U(r) : r = |x-x_0|, x \in \partial D_*\} \mbox{ and } \rho_{*2} =|x_{*2}-x_0|. \label{minimum point}
\end{eqnarray}
 Notice that $\nu(x_{*i})$  equals $\frac{x_{*i}-x_0}{\rho_{*i}}$ for $i=1, 2$.  Also, the case where $\rho_{*1}=\rho_{*2}$ may occur for instance if  $D_*$ is a ball centered at $x_0$. For $r > 0$, we set 
 \begin{equation}
 \label{key function 2}
 \hat{U}(r) 
 = \left\{\begin{array}{rll}
 U(\rho_{*2}) + \frac{\sigma_s}{\sigma_c}\left(U(r)- U(\rho_{*2}) \right)& \mbox{ if } \sigma_s > \sigma_c\ ,
 \\
 U(\rho_{*1}) + \frac{\sigma_s}{\sigma_c}\left(U(r)- U(\rho_{*1})\right) & \mbox{ if } \sigma_s < \sigma_c\ .
\end{array}\right.
\end{equation}
Then, as in \eqref{hat against the original},  it follows that
\begin{equation}
\label{hat is greater than the original}
\hat{U} \ge U \ \mbox{ on } \partial D_*.
\end{equation}
Observe that  
\begin{equation}
\label{hat function is good 2}
-\Delta \hat{U} = \frac 1{\sigma_c}\  \mbox{ and } \ \frac {\partial\hat{U}}{\partial r}= \frac {\sigma_s}{\sigma_c} \frac {\partial U}{\partial r}\ \mbox{ in } \overline{D_*}\setminus\{x_0\},\ \mbox{ and }\ \lim_{x \to x_0}\hat{U} = + \infty.
\end{equation}
Therefore, since we also have \eqref{the equations V satisfies}, it follows from \eqref{hat is greater than the original} and the strong comparison principle that
\begin{equation}
\label{hat is greater than V}
\hat{U} > V \ \mbox{ in }  D_*\setminus \{ x_0\}.
\end{equation}
The transmission condition \eqref{transmission condition} with the definition of $\hat{U}$ tells us that 
\begin{equation}
\label{transmission at the point 2}
\hat{U} = V\ \mbox{ and } \ \frac {\partial \hat{U}}{\partial\nu} = \frac {\partial V}{\partial\nu} \ \mbox{ at } x = x_{*i} \in \partial D_*,
\end{equation}
since $\nu(x_{*i})$ equals $\frac{x_{*i}-x_0}{\rho_{*i}}$ for $i=1, 2$.  Therefore applying  Hopf's boundary point lemma to the harmonic function $\hat{U} - V$ gives a contradiction to \eqref{transmission at the point 2},  and hence case (ii) never occurs.

It remains to show that case (iii) does not occur. In case (iii), since $c_1 < 0$, there exists a unique critical point $r = \rho_c$ of  $U(r)$ such that
\begin{eqnarray}
&&U(\rho_c) = \max\{ U(r) : r > 0 \} > 0\ \mbox{ and } 0 < \rho_c < \rho\ ;  \label{the maximum of U}
\\
&& U^\prime(r) < 0 \ \mbox{ if } r > \rho_c\ \mbox{ and }  U^\prime(r) > 0 \ \mbox{ if } 0< r < \rho_c\ ;  \label{shape of U}
\\
&& \  \lim_{r \to 0} U(r) = - \infty \ \mbox{ and } x_0 \in D. \label{singularity}
\end{eqnarray}
Let us choose the connected component $D_*$ of $D$ satisfying $x_0 \in D_*$.  Then, since $\overline{D_*} \subset \Omega = B_{\rho}(x_0)$,  as in case (ii), we see that there exist $\rho_{*1},  \rho_{*2} \in  (0, \rho)$  and $x_{*1}, x_{*2} \in \partial D_*$ which  satisfy \eqref{maximum point} and \eqref{minimum point}. In view of the shape of the graph of $U$, we have from the transmission condition \eqref{transmission condition} that at $x_{*i} \in \partial D_*, i=1, 2, $
\begin{equation}
\label{vanish normal derivative}
\frac {\partial V}{\partial\nu} = \frac{\sigma_s}{\sigma_c} \frac {\partial U}{\partial\nu}
 = \left\{\begin{array}{rll}
 0 \ &\mbox{ if } \rho_{*i} = \rho_c\ ,
 \\
\frac{\sigma_s}{\sigma_c} U^\prime \ &\mbox{ if } \rho_{*i} \not= \rho_c\ ,
\end{array}\right.
\end{equation}
where, in order to see that  $\nu(x_{*i})$  equals $\frac{x_{*i}-x_0}{\rho_{*i}}$ if $\rho_{*i} \not= \rho_c $, we used the fact that both $D_*$ and $B_{\rho}(x_0) \setminus \overline{D_*}$ are connected and $x_0 \in D_*$.  Also, the case where $\rho_{*1}=\rho_{*2}$ may occur for instance if  $D_*$ is a ball centered at $x_0$. For $r > 0$, we define $\hat{U} = \hat{U}(r)$ by 
 \begin{equation}
 \label{key function 3}
 \hat{U}(r) 
 = \left\{\begin{array}{rll}
 U(\rho_{*1}) + \frac{\sigma_s}{\sigma_c}\left(U(r)- U(\rho_{*1}) \right)& \mbox{ if } \sigma_s > \sigma_c\ ,
 \\
 U(\rho_{*2}) + \frac{\sigma_s}{\sigma_c}\left(U(r)- U(\rho_{*2})\right) & \mbox{ if } \sigma_s < \sigma_c\ .
\end{array}\right.
\end{equation}
Remark that \eqref{key function 3} is opposite to \eqref{key function 2}. Then, as in \eqref{hat is greater than the original},  it follows that
\begin{equation}
\label{hat is less than the original}
\hat{U} \le U \ \mbox{ on } \partial D_*.
\end{equation}
Hence, by proceeding with the strong comparison principle as in case (ii), we conclude that
\begin{equation}
\label{hat is less than V}
\hat{U} < V \ \mbox{ in }  D_*\setminus \{ x_0\}.
\end{equation}
Then, it follows from the definition of $\hat{U}$ and \eqref{vanish normal derivative} that \eqref{transmission at the point 2} also holds true. In conclusion,  applying  Hopf's boundary point lemma to the harmonic function $\hat{U} - V$ gives a contradiction to \eqref{transmission at the point 2},  and hence case (iii) never occurs.

\setcounter{equation}{0}
\setcounter{theorem}{0}

\section{Proof of Theorem \ref{th:stationary isothermic cauchy}}
\label{section4}

Let $u$ be the solution of  problem  \eqref{heat Cauchy}  for $N \ge 3$. For assertion (2) of Theorem \ref{th:stationary isothermic cauchy},
with the aid of Aleksandrov's sphere theorem \cite[p. 412]{Alek1958vestnik},  Lemma \ref{le: constant weingarten curvature} yields that $\gamma$ and $\Gamma$ are concentric spheres. Denote by $x_0 \in \mathbb R^N$ the common center of $\gamma$ and $\Gamma$. By combining the initial  condition of problem  \eqref{heat Cauchy}  and  the assumption \eqref{stationary isothermic surface partially} with the real analyticity in $x$ of $u$ over $\mathbb R^N \setminus\overline{D}$ coming from $\sigma_s=\sigma_m$,  we see that $u$ is radially symmetric with respect to $x_0$  in $x$ on $\left(\mathbb R^N\setminus\overline{D}\right) \times (0,\infty)$. Here we used the assumption that $\Omega\setminus\overline{D}$ is connected. Moreover, in view of the initial  condition of problem  \eqref{heat Cauchy},  we can distinguish the following two cases as in section \ref{section3}:
$$
\mbox{\rm (I)  } \Omega \mbox{ is a ball;}\qquad \mbox{\rm (II)  } \Omega \mbox{ is a spherical shell.}
$$ 
For assertion (1) of Theorem \ref{th:stationary isothermic cauchy}, with the aid of Aleksandrov's sphere theorem \cite[p. 412]{Alek1958vestnik},  Lemma \ref{le: constant weingarten curvature 2} yields that $\partial G$ and $\partial\Omega$ are concentric spheres, since every component of $\partial\Omega$ is a sphere with the same curvature. Therefore, only the case (I) remains for assertion (1) of Theorem \ref{th:stationary isothermic cauchy}. Also, denoting by $x_0 \in \mathbb R^N$ the common center of $\partial G$ and $\partial\Omega$ and combining the initial  condition of problem  \eqref{heat Cauchy}  and  the assumption \eqref{stationary isothermic surface} with the real analyticity in $x$ of $u$ over $\Omega \setminus\overline{D}$ yield that  $u$ is radially symmetric with respect to $x_0$  in $x$ on $\left(\mathbb R^N\setminus\overline{D}\right) \times (0,\infty)$. 

By virtue of (2) of Lemma \ref{le:initial behavior exponential decay and decay at infinity},  since $N \ge 3$,  we can introduce the following three auxiliary functions $U = U(x), \ V = V(x)$ and $W =W(x)$ by
\begin{eqnarray}
&& U(x) = \int_0^\infty(1-u(x,t) )\ dt\quad \mbox{ for } x \in \Omega \setminus\overline{D}, \label{auxiliary function on shell 2}
\\
&& V(x) = \int_0^\infty(1-u(x,t) )\ dt\quad \mbox{ for } x \in D,\label{auxiliary function on core 2}
\\
&& W(x) = \int_0^\infty(1-u(x,t) )\ dt\quad \mbox{ for } x \in \mathbb R^N \setminus \overline{\Omega}.\label{auxiliary function on medium}
\end{eqnarray}
Then we observe that
\begin{eqnarray}
&&-\Delta U = \frac 1{\sigma_s} \ \mbox{ in } \Omega \setminus\overline{D}, \ -\Delta V = \frac 1{\sigma_c} \ \mbox{ in } D,
 \ -\Delta W = 0\ \mbox{ in } \mathbb R^N\setminus\overline{\Omega}, \label{poisson and Laplace equations}
\\
&&U = V\ \mbox{ and }\  \sigma_s \frac {\partial U}{\partial\nu} = \sigma_c \frac {\partial V}{\partial\nu}  \ \mbox{ on } \partial D, \label{transmission condition between U and V}
\\
&&U = W\ \mbox{ and }\  \sigma_s \frac {\partial U}{\partial\nu} = \sigma_m \frac {\partial W}{\partial\nu}  \ \mbox{ on } \partial\Omega, \label{transmission condition between U and W}
\\
&& \lim_{|x| \to \infty} W(x) = 0, \label{decay at infinity}
\end{eqnarray}
where $\nu = \nu(x)$ denotes the outward unit normal vector to $\partial D$ at $x \in \partial D$ or to $\partial\Omega$ at $x \in \partial\Omega$ and \eqref{transmission condition between U and V} - \eqref{transmission condition between U and W} are the transmission conditions. Here we used (4) of Lemma \ref{le:initial behavior exponential decay and decay at infinity} to obtain \eqref{decay at infinity}.

Let us follow the proof of Theorem \ref{th:stationary isothermic}. We first show that case (II)  for  assertion (2) of Theorem \ref{th:stationary isothermic cauchy} does not occur. Set $\Omega = B_{\rho_+}(x_0) \setminus \overline{B_{\rho_-}(x_0)}$ for some numbers $\rho_+ > \rho_- > 0$. Since $u$ is radially symmetric with respect to $x_0$  in $x$ on $\left(\mathbb R^N\setminus\overline{D}\right) \times (0,\infty)$, we can obtain from \eqref{poisson and Laplace equations}-\eqref{decay at infinity} that for $r =|x-x_0| \ge 0$
\begin{eqnarray*}
& U = c_1 r^{2-N} - \frac 1{2N\sigma_s} r^2 + c_2\ &\mbox{ for } \rho_-\le r \le \rho_+,
\\
& W = c_3 r^{2-N} \ &\mbox{ for } r \ge \rho_+,
\\
& W = c_4 \ &\mbox{ for } 0 \le r \le \rho_-,
\end{eqnarray*}
where $c_1, \dots, c_4$ are some constants, since $\Omega \setminus \overline{D}$ is connected. Remark that $U$ can be extended as a radially symmetric function of $r$ in $\mathbb R^N \setminus \{x_0\}$.
We observe that $c_4 > 0$ and $c_3 > 0$. Also it follows from \eqref{transmission condition between U and W} that
$U^\prime(\rho_-) = 0$ and $U^\prime(\rho_+) < 0$, and hence 
$$
c_1 < 0\ \mbox{ and } \ U^\prime < 0 \mbox{ on } (\rho_-,\rho_+].
$$
Then the same argument as in the corresponding case in the proof of Theorem \ref{th:stationary isothermic} works and a contradiction to the transmission condition \eqref{transmission condition between U and V} can be obtained. Thus case (II)  for  assertion (2) of Theorem \ref{th:stationary isothermic cauchy} never occurs.

Let us proceed to case (I). Set $\Omega = B_\rho(x_0)$ for some number $\rho > 0$. Since $u$ is radially symmetric with respect to $x_0$  in $x$ on $\left(\mathbb R^N\setminus\overline{D}\right) \times (0,\infty)$, we can obtain from \eqref{poisson and Laplace equations}-\eqref{decay at infinity} that for $r =|x-x_0| \ge 0$
\begin{eqnarray*}
& U = c_1 r^{2-N} - \frac 1{2N\sigma_s} r^2 + c_2\ &\mbox{ for } x \in \overline{\Omega} \setminus D,
\\
& W = c_3 r^{2-N} \ &\mbox{ for } r \ge \rho,
\end{eqnarray*}
where $c_1, c_2, c_3$ are some constants, since $\Omega \setminus \overline{D}$ is connected. Remark that $U$ can be extended as a radially symmetric function of $r$ in $\mathbb R^N \setminus \{x_0\}$. Therefore it follows from \eqref{transmission condition between U and W} that $U^\prime(\rho) < 0$. As in the proof of Theorem \ref{th:stationary isothermic}, 
We distinguish the following three cases:
$$
\mbox{\rm (i)  } c_1 = 0;\qquad \mbox{\rm (ii)  } c_1 > 0; \qquad \mbox{\rm (iii)  } c_1 < 0.
$$ 
Because of the fact that $U^\prime(\rho) < 0$, the same arguments as  in the proof of Theorem \ref{th:stationary isothermic} works to conclude that only case (i) occurs and $D$ must be a ball centered at $x_0$. \qed

\vskip 4ex
\bigskip
\noindent{\large\bf Acknowledgement.}
\smallskip

The main results of the present paper were discovered while the author was visiting 
the National Center for Theoretical Sciences (NCTS) Mathematics Division in the National
Tsing Hua University; he wishes to thank NCTS for its kind hospitality.


\end{document}